# An integrated approach to grey relational analysis, analytic hierarchy process and data envelopment analysis

Mohammad Sadegh Pakkar
*Department of Finance and Operations, Laurentian University, Ontario, Canada*

## Abstract

**Purpose** – This paper aims to propose an integration of the analytic hierarchy process (AHP) and data envelopment analysis (DEA) methods in a multiattribute grey relational analysis (GRA) methodology in which the attribute weights are completely unknown and the attribute values take the form of fuzzy numbers.

**Design/methodology/approach** – This research has been organized to proceed along the following steps: computing the grey relational coefficients for alternatives with respect to each attribute using a fuzzy GRA methodology. Grey relational coefficients provide the required (output) data for additive DEA models; computing the priority weights of attributes using the AHP method to impose weight bounds on attribute weights in additive DEA models; computing grey relational grades using a pair of additive DEA models to assess the performance of each alternative from the optimistic and pessimistic perspectives; and combining the optimistic and pessimistic grey relational grades using a compromise grade to assess the overall performance of each alternative.

**Findings** – The proposed approach provides a more reasonable and encompassing measure of performance, based on which the overall ranking position of alternatives is obtained. An illustrated example of a nuclear waste dump site selection is used to highlight the usefulness of the proposed approach.

**Originality/value** – This research is a step forward to overcome the current shortcomings in the weighting schemes of attributes in a fuzzy multiattribute GRA methodology.

**Keywords** Data envelopment analysis, Grey relational analysis, Analytic hierarchy process, Multiple attribute decision making, Fuzzy numbers, Weighting

## Introduction

Grey relational analysis (GRA) is a part of the grey system theory proposed by Deng (1982), which is suitable for solving a variety of multiple attribute decision making

**JEL classification** – C02, C61, D81, M11

The author expresses his gratitude to the Editor-in-Chief and the anonymous reviewers whose valuable comments helped to improve the article.

(MADM) problems with both crisp and fuzzy data (Goyal and Grover, 2012; Wei *et al.*, 2011; Wei, 2010; Hou, 2010; Olson and Wu, 2006). Grey sets can be considered as an extension to fuzzy sets by restricting the characteristic (or membership) function values of a set within [0,1] (Yang and John, 2012; Li *et al.*, 2008). GRA solves MADM problems by aggregating multiple attribute values, which are usually incommensurable, into a single value for each alternative (Kuo *et al.*, 2008). In the traditional GRA method, assigning equal weights to attributes for each alternative is a norm. Nevertheless, the validity of using the equal weight assumption for all the alternatives to be assessed can be questioned, as each one of these has its own characteristics and preferences. Therefore, the study on attribute weighting can be an interesting, but it is a controversial topic in the field of GRA. Fortunately, the development of modern operational research has provided us with two excellent tools, namely, analytic hierarchy process (AHP) and data envelopment analysis (DEA), which can be used to derive attribute weights for use in GRA.

AHP is a subjective data-oriented procedure, which can reflect the relative importance of a set of attributes and alternatives based on the formal expression of the decision-maker's preferences. AHP usually involves three basic functions: structuring complexities, measuring on a ratio-scale and synthesizing (Saaty, 1987). Some researchers incorporate fuzzy set theory in the conventional AHP to express the uncertain comparison judgments as fuzzy numbers (Osman *et al.*, 2013; Javanbarg *et al.*, 2012; Kahraman *et al.*, 2003). However, AHP has been criticized because of the arbitrary nature of the ranking process (Ahmad *et al.*, 2006; Swim, 2001; Dyer *et al.*, 1990). In fact, the AHP weights are based on the experts' personal experiences and their subjective judgments. If the selection of experts is different, then the weights obtained will be different (Liu, 2009; Liu and Chen, 2004). The application of AHP with GRA can be seen in the studies by Birgün and Güngör (2014), Jia *et al.* (2011) and Zeng *et al.* (2007).

Alternatively, DEA is an objective data-oriented approach to assess the relative performance of a group of decision-making units (DMUs) with multiple inputs and outputs (Cooper *et al.*, 2011). Traditional DEA models require crisp input and output data. However, in recent years, fuzzy set theory has been proposed for quantifying imprecise and vague data in DEA models (Hatami-Marbini *et al.*, 2013; Wen and Li, 2009; Lertworasirikul *et al.*, 2003). In the field of GRA, DEA models without explicit inputs are applied, that is the models in which only pure outputs or index data are taken into account (Liu *et al.*, 2011). The other combined GRA and DEA methodologies can be found in the literature, such as using GRA for the selection of inputs and outputs in DEA (Wang *et al.*, 2010; Bruce Ho, 2011), using GRA for ranking efficient DMUs in DEA with crisp data (Girginer *et al.*, 2015), using GRA for ranking DMUs in DEA with grey data, that is the unknown numbers, which have clear upper and lower limits (Markabi and Sabbagh, 2014) and weighting GRA using a cross-efficiency model (Markabi and Sarbijan, 2015). However, the main problem of using traditional DEA models in GRA is that several alternatives may receive a grey relational grade of 1, which means that all of these alternatives are ranked in the first position (Jun and Xiaofei, 2013; Wu and Olson, 2010). To overcome this deficiency, Zheng and Lianguang (2013) propose a super-efficiency model, based on the super-efficiency ranking method of Andersen and Petersen (1993). According to this model, these alternatives are allowed to obtain a grey relational grade greater than 1 by removing the constraint that bounds the grade of the assessed alternative. These grades are then used to rank the alternatives and thereby eliminate some of the ties that occur in the selection of the best alternative. As shown in all the aforementioned models, each alternative has been allowed to choose its own most favorable weights to maximize its performance, they can be called optimistic DEA models.

On the other hand, a similar approach can be developed to assess the performance of each alternative under the least favorable weights, which, in fact, is pessimistic. According to this approach, each alternative is compared with the worst alternatives and is assessed by its grey relational grade as the ratio of the distance from the worst frontier. It is worth pointing out that the worst-practice frontier approach is not a new approach in the DEA literature. Conceptually, it is parallel to the worst possible efficiency concept as discussed by Wang and Luo (2006), Takamura and Tone (2003), Jahanshahloo and Afzalinejad (2006) and Liu and Chen (2009). Nevertheless, as far as we know, this approach has never been applied to the field of GRA.

It can be argued that both optimistic and pessimistic DEA approaches should be considered together to obtain attribute weights in GRA, and any approach considers that only one of them is biased (Wang *et al.*, 2007). In both approaches, each DMU or alternative can freely choose its own system of weights to optimize its performance. However, this freedom of choosing weights is equivalent to keeping the preferences of a decision-maker out of the decision process. In fact, an alternative may be indicated as the best (worst) one by assigning zero values to the weights of some attributes and neglecting the relative priorities of these attributes in the decision-making process.

To overcome these issues, we propose the integration of AHP and DEA to obtain the attribute weights in GRA from both the optimistic and pessimistic perspectives. This can be implemented by imposing weight restrictions in DEA-based GRA models, using AHP (Pakkar, 2016a). There are two types of weight restrictions: homogeneous and non-homogeneous (Podinovski, 2004). A homogeneous weight restriction has a zero-free constant on its right-hand side (RHS), while a non-homogeneous weight restriction has a non-zero constant on its RHS. The bounds on the relative values of input and output weights (Lee *et al.*, 2012; Liu *et al.*, 2005; Tseng *et al.*, 2009; Kong and Fu, 2012) and those on the virtual weights of inputs and outputs (Premachandra, 2001; Shang and Sueyoshi, 1995) are typical examples of homogeneous weight restrictions using AHP. The absolute weight bounds that directly impose bounds on input and output weights (Entani *et al.*, 2004; Pakkar, 2016b) belong to the type of non-homogeneous weight restrictions using AHP. The DEA models discussed in this paper are all based on the additive models (Cooper *et al.*, 1999) without explicit inputs in which output data are defined by a fuzzy GRA method. In these models, the decision makers are allowed to reflect those considerations not embodied in data. This would allow us to easily incorporate AHP weights into the multiplier form of additive models using suitable non-homogeneous weight restrictions.

## Methodology

As has been mentioned earlier, a critical issue in using the GRA method is the subjectivity in assigning weights for attributes. As different weight combinations may lead to different ranking results, it is unlikely that all the decision-makers would easily reach a consensus in determining an appropriate set of weights. In addition, it may not be easy to obtain expert information for deriving the weights. Although the use of equal weights seems to be a relatively fair choice, some decision-makers may still have different opinions on the relative importance of attributes. To avoid these issues, an integration of AHP and DEA models is given here to obtain the attribute weights in a fuzzy GRA methodology. This can be implemented through the following steps (Figure 1):

- computing the grey relational coefficients for alternatives with respect to each attribute using a fuzzy GRA method. Grey relational coefficients provide the required (output) data for additive DEA models;
- computing the priority weights of attributes using the AHP method to impose weight bounds on attribute weights in additive DEA models;
- computing grey relational grades using a pair of additive DEA models to assess the performance of each alternative from the optimistic and pessimistic perspectives; and
- combining the optimistic and pessimistic grey relational grades using a compromise grade to assess the overall performance of each alternative.

*Fuzzy grey relational analysis*

GRA can be applied to both crisp and fuzzy data. Here, we use it as a means to obtain a solution from fuzzy data. Let $A = \{A_1, A_2, \ldots, A_m\}$ be a discrete set of alternatives and $C = \{C_1, C_2, \ldots, C_n\}$ be a set of attributes. Let $\tilde{y}_{ij} = (y_{1ij}, y_{2ij}, y_{3ij}, y_{4ij})$ be a trapezoidal fuzzy number representing the value of attribute $C_j (j = 1, 2, \ldots, n)$ for alternative $A_i\,(i = 1, 2, \ldots, m)$. Using the $\alpha$-cut technique, a trapezoidal fuzzy number can be transformed into an interval number as follows:

$$y_{ij} = [y_{ij}^-, y_{ij}^+] = [\alpha y_{2ij} + (1 - \alpha) y_{1ij}, \alpha y_{3ij} + (1 - \alpha) y_{4ij})] \quad (1)$$

where $y_{ij} = [y_{ij}^-, y_{ij}^+], y_{ij}^- \leq y_{ij}^+$, is an interval number representing the value of attribute $C_j (j = 1, 2, \ldots, n)$ for alternative $A_i\,(i = 1, 2, \ldots, m)$. Then alternative $A_i$ is characterized by a vector $Y_i = ([y_{i1}^-, y_{i1}^+], [y_{i2}^-, y_{i2}^+], \dots, [y_{in}^-, y_{in}^+])$ of attribute values. The term $Y_i$ can be translated into the comparability sequence $R_i = ([r_{i1}^-, r_{i1}^+], [r_{i2}^-, r_{i2}^+], \dots, [r_{in}^-, r_{in}^+])$ by using the following equations (Zhang *et al.*, 2005):

$$[r_{ij}^-, r_{ij}^+] = \left[\frac{y_{ij}^-}{y_{j(\max)}^+}, \frac{y_{ij}^+}{y_{j(\max)}^+}\right] \forall j,\ y_{j(\max)}^+ = \max\{y_{1j}^+, y_{2j}^+, \dots, y_{mj}^+\} \text{ for desirable attributes,} \tag{2}$$

$$[r_{ij}^-, r_{ij}^+] = \left[\frac{y_{i(\min)}^-}{y_{ij}^-}, \frac{y_{i(\min)}^-}{y_{ij}^+}\right] \forall j,\ y_{j(\min)}^- = \min\{y_{1j}^-, y_{2j}^-, \dots, y_{mj}^-\} \text{ for undesirable attributes.} \tag{3}$$

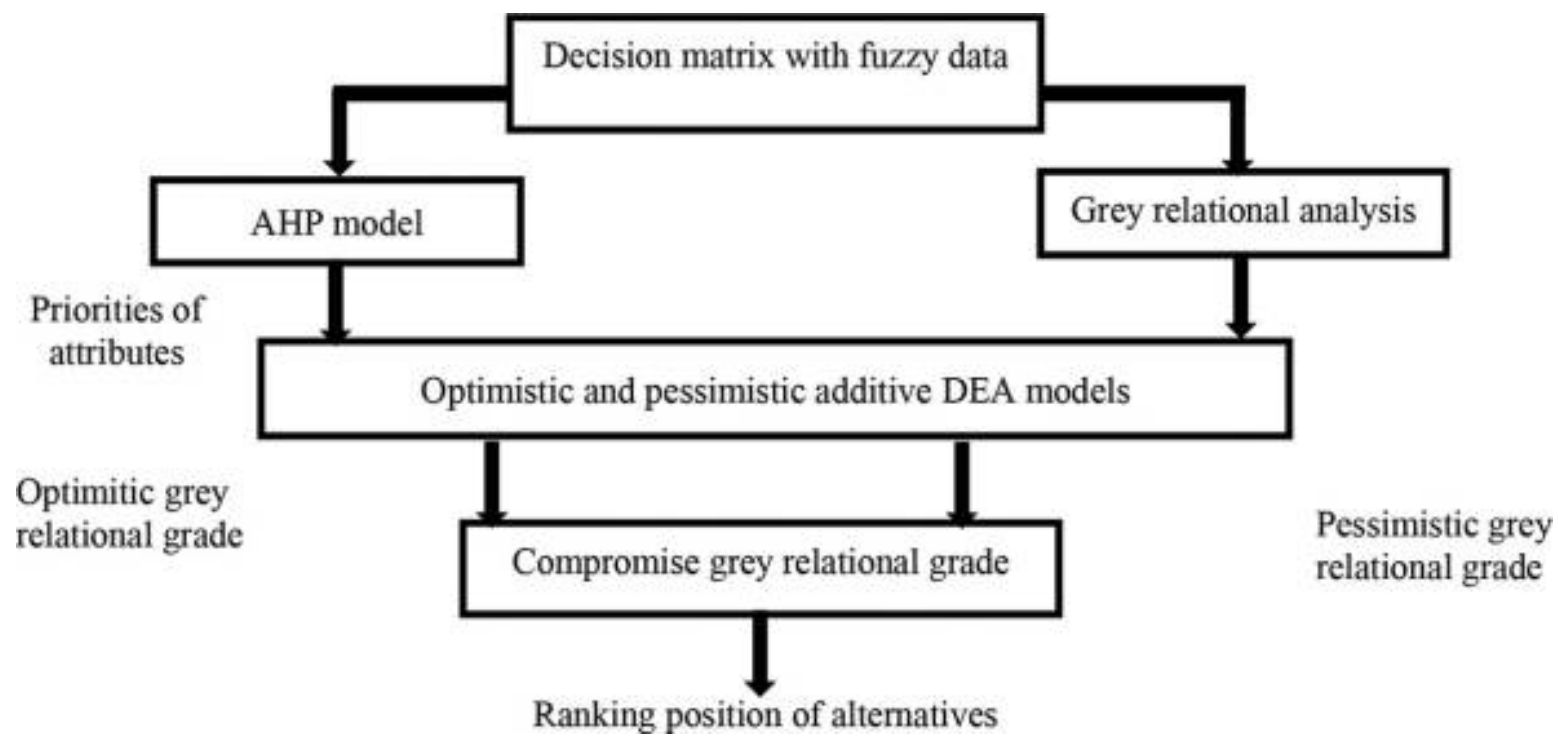

**Figure 1.**
An integrated approach to GRA, AHP and DEA

Now, let $A_0$ be a virtual ideal alternative, which is characterized by a reference sequence $U_0 = ([u_{01}^-, u_{01}^+], [u_{02}^-, u_{i2}^+], \dots, [u_{0n}^-, u_{0n}^+])$ of the maximum attribute values as follows:

$$u_{0j}^- = \max\{r_{1j}^-, r_{2j}^-, \dots, r_{mj}^-\}\ \forall j, \tag{4}$$

$$u_{0j}^+ = \max\{r_{1j}^+, r_{2j}^+, \dots, r_{mj}^+\}\ \forall j. \tag{5}$$

To measure the degree of similarity between $r_{ij} = [r_{ij}^-, r_{ij}^+]$ and $u_{0j} = [u_{0j}^-, u_{0j}^+]$ for each attribute, the grey relational coefficient, $\xi_{ij}$, can be calculated as follows:

$$\xi_{ij} = \frac{\min_i \min_j |[u_{0j}^-, u_{0j}^+] - [r_{ij}^-, r_{ij}^+]| + \rho \max_i \max_j |[u_{0j}^-, u_{0j}^+] - [r_{ij}^-, r_{ij}^+]|}{|[u_{0j}^-, u_{0j}^+] - [r_{ij}^-, r_{ij}^+]| + \rho \max_i \max_j |[u_{0j}^-, u_{0j}^+] - [r_{ij}^-, r_{ij}^+]|} \tag{6}$$

while the distance between $u_{0j} = [u_{0j}^-, u_{0j}^+]$ and $r_{ij} = [r_{ij}^-, r_{ij}^+]$ is measured by $|u_{0j} - r_{ij}| = \max(|u_{0j}^- - r_{ij}^-|, |u_{0j}^+ - r_{ij}^+|)$. $\rho \in [0, 1]$ is the distinguishing coefficient, which is generally $\rho = 0.5$. It should be noted that the final results of GRA for MADM problems are very robust to changes in the values of $\rho$. Therefore, selecting the different values of $\rho$ would only slightly change the rank order of attributes (Kuo *et al.*, 2008). To find an aggregated measure of similarity between alternative $A_i$, characterized by the comparability sequence $R_i$, and the ideal alternative $A_0$, characterized by the reference sequence $U_0$, over all the attributes, the grey relational grade, $\Gamma_i$, can be computed as follows:

$$\Gamma_i = \sum_{j=1}^{n} w_j \xi_{ij} \tag{7}$$

where $w_j$ is the weight of attribute $C_j$. In the next section, we show how the AHP model can be used to obtain the priority weights of attributes for each alternative.

*The analytic hierarchy process*

The AHP procedure for computing the priority weights of attributes may be broken down into the following steps:

*Step 1*: A decision-maker makes a pairwise comparison matrix of different attributes, denoted by $B$ with the entries of $b_{hq}$ $(h = q = 1, 2, \ldots, n)$. The comparative importance of attributes is provided by the decision-maker using a rating scale. Saaty (1987) recommends using a 1-9 scale.

*Step 2*: The AHP method obtains the priority weights of attributes by computing the eigenvector of matrix $B$ [equation (8)], $e = (e_1, e_2, \ldots, e_j)^T$, which is related to the largest eigenvalue, $\gamma_{\max}$:

$$Be = \gamma_{\max} e \tag{8}$$

To determine whether the inconsistency in a comparison matrix is reasonable the random consistency ratio, $C.R.$, can be computed by the following equation:

$$C.R. = \frac{\gamma_{\max} - N}{(N - 1)R.I.} \tag{9}$$

where $R.I.$ is the average random consistency index and $N$ is the size of a comparison matrix.

### *Optimistic and pessimistic additive DEA models*

As all the grey relational coefficients are benefit (output) data, an optimistic additive DEA model for obtaining attribute weights in GRA can be developed like the additive model in Cooper *et al.* (1999) without explicit inputs as follows:

$$\begin{aligned} P_k &= \max \sum_{j=1}^{n} e_j s_j^+ \\ \text{s.t.} & \sum_{i=1}^{m} \lambda_i \xi_{ij} - s_j^+ = \xi_{kj} \ \forall j, \\ & \sum_{i=1}^{m} \lambda_i = 1 \\ & s_j^+, \lambda_i \geq 0, \end{aligned} \tag{10}$$

where $1 - P_k$ indicates the grey relational grade, $\Gamma_k (k = 1, 2, \ldots, m)$, for alternative under assessment $A_k$ (known as a DMU in the DEA terminology) and $0 \leq P_k \leq 1$. $s_j^+$ is the slack variable of attribute $C_j (j = 1, 2, \ldots, n)$, expressing the difference between the performance of a composite alternative and the performance of the assessed alternative with respect to each attribute. In other words, $s_j^+$ identifies a shortfall in the attribute value of $C_j$ for alternative $A_k$. Obviously, when $P_k = 0$ alternative $A_k$ is considered as the best alternative in comparison to all the other alternatives. $e_j$ is the priority weight of attribute $C_j$, which is defined out of the internal mechanism of DEA using AHP and $\lambda_i$ is the weight of alternative $A_i (i = 1, 2, \ldots, m)$. The convexity constraint in model (10) meets the assumption of *variable returns-to-scale* (VRS) frontier for an additive model. The dual of model (10) can be developed as follows:

$$
\begin{aligned}
\Gamma_k &= \max \sum_{j=1}^{n} w_j \xi_{kj} - w_0 \\
\text{s.t.} & \sum_{j=1}^{n} w_j \xi_{ij} - w_0 \leq 1 \; \forall i, \\
& w_j \geq e_j \; \forall j, \\
& w_0 \; free.
\end{aligned} \tag{11}
$$

This model is useful for our purpose in dealing with grey relational grades. The objective function in model (11) maximizes the ratio of the grey relational grade of alternative $A_k$ to the maximum grey relational grade across all alternatives for the same set of weights ($\max \Gamma_k / \max_{i=1,\dots,m} \Gamma_i$), while the priority weights obtained by AHP impose the lower bounds on the attribute weights. Hence, an optimal set of weights in model (11) represents $A_k$ in the best light compared to all the other alternatives while it reflects a *priori* information about the priorities of attributes, simultaneously. Finally, one should notice that the optimistic additive DEA models bounded by AHP does not necessarily yield results that are different from those obtained from the original additive DEA models (Charnes *et al.*, 1985). In particular, it does not increase the power of discrimination between the considerable number of alternatives, which are usually ranked in the first place by obtaining the grey relational grades of 1. To overcome these issues, we develop the additive models from the pessimistic point of view in which each alternative is assessed based on its distance from the worst practice frontier as follows:

$$
\begin{aligned}
P'_k &= \max \sum_{j=1}^{n} e_j s_j^- \\
\text{s.t.} & \sum_{i=1}^{m} \lambda'_i \xi_{ij} + s_j^- = \xi_{kj} \; \forall j, \\
& \sum_{i=1}^{m} \lambda'_i = 1 \\
& s_j^-, \lambda'_i \geq 0,
\end{aligned} \tag{12}
$$

Note that the only difference between model (10) and model (12) is the signs of slack variables in the first set of constraint. $s_j^-$ is the slack variable of attribute $C_j (j = 1, 2, \dots, n)$, expressing the difference between the performance of the assessed alternative and the performance of a composite alternative with respect to each attribute. In other words, $s_j^-$ identifies the excess values of attribute $C_j$ for alternative $A_k$. This is obvious when $P'_k = 0$ alternative $A_k$ is considered as the worst alternative compared to all the other alternatives. In some cases, the worst alternatives, however, may also be the best alternatives. This happens when the assessed alternative is the best in some attributes, while it is the worst in some other attributes. The dual of model (12) can be developed as follows:

$$\Gamma_k^{'} = \min \sum_{j=1}^{n} w_j^{'} \xi_{kj} + w_0^{'}$$
$$\text{s.t.} \sum_{j=1}^{n} w_j^{'} \xi_{ij} + w_0^{'} \geq 1 \ \forall i,$$
$$w_j^{'} \geq e_j \forall j,$$
$$w_0^{'} \ free. \tag{13}$$

Here, we seek the worst weights in the sense that the objective function in model (13) is minimized. The first set of constraints assures that the computed weights do not attain a grade smaller than 1. Each alternative is compared with these worst alternatives and is assessed based on the ratio of the distance from the worst-practice frontier. It is worth pointing out that the pessimistic additive models in this paper are not brand-new models in the DEA literature. Conceptually, it is parallel to the additive DEA models as discussed by Jahanshahloo and Afzalinejad (2006) for ranking alternatives on a full inefficient-frontier. Nevertheless, as far as we know, it is the first time that they are applied to the field of GRA.

To combine the grey relational grades obtained from models (11) and (13), that is the best and worst sets of weights, the linear combination of corresponding normalized grades is recommended as follows (Zhou *et al.*, 2007):

$$\Delta_k(\beta) = \beta \frac{\Gamma_k - \Gamma_{\min}}{\Gamma_{\max} - \Gamma_{\min}} + (1 - \beta) \frac{\Gamma_k^{'} - \Gamma_{\min}^{'}}{\Gamma_{\max}^{'} - \Gamma_{\min}^{'}} \tag{14}$$

where $\Gamma_{\max} = \max\{\Gamma_k, k = 1, 2, \ldots, m\}$, $\Gamma_{\min} = \min\{\Gamma_k, k = 1, 2, \ldots, m\}$, $\Gamma_{\max}^{'} = \max\{\Gamma_k^{'}, k = 1, 2, \ldots, m\}$, $\Gamma_{\min}^{'} = \min\{\Gamma_k^{'}, k = 1, 2, \ldots, m\}$ and $0 \leq \beta \leq 1$ is an adjusting parameter, which may reflect the preference of a decision-maker on the best and worst sets of weights. $\Delta_k(\beta)$ is a normalized compromise grade in the range [0,1].

*Numerical example: nuclear waste dump site selection*

In this section, we present the application of the proposed approach for nuclear waste dump site selection. The multi-attribute data, adopted from Wu and Olson (2010), are presented in Table I. There are 12 alternative sites and 4 performance attributes. Cost, lives lost and risk are undesirable attributes and civic improvement is a desirable attribute. Cost is in billions of dollars. Lives lost reflect expected lives lost from all exposures. Risk shows the risk of catastrophe (earthquake, flood, etc.) and civic improvement is the improvement of the local community due to the construction and operation of each site. Cost and lives lost are crisp values as outlined in Table I, but risk and civic improvement have fuzzy data for each nuclear dump site.

We use the processed data as reported by Wu and Olson (2010). First the trapezoidal fuzzy data are used to express linguistic data in Table I. Using the $\alpha$-cut technique, the raw data are expressed in fuzzy intervals as shown in Table II. These data are turned into the comparability sequence by using equations (2) and (3). Each attribute is now on a common 0-1 scale where 0 represents the worst imaginable attainment on an attribute, and 1 represents the best possible attainment.

Table III shows the results of a pairwise comparison matrix in the AHP model as constructed by the author in Expert Choice software. The priority weight for each attribute would be the average of the elements in the corresponding row of the normalized matrix of pairwise comparison, shown in the last column of Table III. One can argue that the priority weights of attributes must be judged by nuclear safety experts. However, as the aim of this section is just to show the application of the proposed approach on numerical data, we see no problem to use our judgment alone.

| Site | Cost | Lives | Risk | Civic |
|---|---|---|---|---|
| Nome | 40 | 60 | Very high | Low |
| Newark | 100 | 140 | Very low | Very high |
| Rock Springs | 60 | 40 | Low | High |
| Duquesne | 60 | 40 | Medium | Medium |
| Gary | 70 | 80 | Low | Very high |
| Yakima | 70 | 80 | High | Medium |
| Turkey | 60 | 70 | High | High |
| Wells | 50 | 30 | Medium | Medium |
| Anaheim | 90 | 130 | Very high | Very low |
| Epcot | 80 | 120 | Very low | Very low |
| Duckwater | 80 | 70 | Medium | Low |
| Santa Cruz | 90 | 100 | Very high | Very low |

**Table I.** Data for nuclear waste dump site selection

| Site | Cost | Lives lost | Risk | Civic |
|---|---|---|---|---|
| Nome | [0.80-1.00] | [0.40-0.70] | [0.00-0.10] | [0.10-0.30] |
| Newark | [0.00-0.05] | [0.00-0.05] | [0.90-1.00] | [0.90-1.00] |
| Rock Springs | [0.70-0.95] | [0.70-0.90] | [0.70-0.90] | [0.70-0.90] |
| Duquesne | [0.50-0.85] | [0.70-0.90] | [0.40-0.60] | [0.40-0.60] |
| Gary | [0.40-0.60] | [0.10-0.30] | [0.70-0.90] | [0.90-1.00] |
| Yakima | [0.50-0.70] | [0.10-0.30] | [0.10-0.30] | [0.40-0.60] |
| Turkey | [0.75-0.90] | [0.20-0.40] | [0.10-0.30] | [0.70-0.90] |
| Wells | [0.85-0.95] | [0.85-1.00] | [0.40-0.60] | [0.40-0.60] |
| Anaheim | [0.00-0.30] | [0.00-0.10] | [0.00-0.10] | [0.00-0.10] |
| Epcot | [0.10-0.40] | [0.00-0.20] | [0.90-1.00] | [0.00-0.10] |
| Duckwater | [0.30-0.50] | [0.20-0.40] | [0.40-0.60] | [0.10-0.30] |
| Santa Cruz | [0.10-0.40] | [0.10-0.30] | [0.00-0.10] | [0.00-0.10] |

**Table II.** Fuzzy interval nuclear waste dump site data

| Attribute | Cost | Lives | Risk | Civic | Priority |
|---|---|---|---|---|---|
| Cost | 1 | 1/5 | 1/2 | 3 | 0.131 |
| Lives | 5 | 1 | 2 | 9 | 0.545 |
| Risk | 2 | 1/2 | 1 | 6 | 0.275 |
| Civic | 1/3 | 1/9 | 1/6 | 1 | 0.05 |

**Note:** $C.R. = 0.01$

**Table III.** Pairwise comparison matrix of four attributes

Using equation (6), all grey relational coefficients are computed to provide the required (output) data for additive DEA models as shown in Table IV. Note that the grey relational coefficients depend on the distinguishing coefficient $\rho$, which here is 0.80.

Table V presents the results obtained from models (11) and (13), as well as the corresponding composite grades at $\beta = 0.5$. If decision-makers have no strong preference, $\beta = 0.5$ would be a fairly neutral and reasonable choice. It can be seen from Table V, the Wells site, with a compromised grade of 1, stands in the first place, while six other alternatives are ranked in the first position by model (11). It is likely because of the fact that the Wells site not only has relatively high values of grey relational coefficients but also has a better combination among the different attributes. This indicates that the proposed approach can significantly improve the degree of discrimination among alternatives. We can also observe that Newark is the best alternative from the optimistic point of view, but it is also the worst alternative from the pessimistic point of view. It is due to the fact that Newark is the best with respect to risk and civic improvement, while it is the worst with respect to

**Table IV.** Results of grey relational coefficients for nuclear waste dump site selection

| Site | Cost | Lives lost | Risk | Civic |
|---|---|---|---|---|
| Nome | 0.9383 | 0.6281 | 0.4578 | 0.4872 |
| Newark | 0.4444 | 0.4444 | 1 | 1 |
| Rock Springs | 0.8352 | 0.8352 | 0.7917 | 0.7917 |
| Duquesne | 0.6847 | 0.8352 | 0.6032 | 0.6032 |
| Gary | 0.6281 | 0.5033 | 0.7917 | 1 |
| Yakima | 0.6847 | 0.5033 | 0.4872 | 0.6032 |
| Turkey | 0.8837 | 0.539 | 0.4872 | 0.7917 |
| Wells | 0.9383 | 1 | 0.6032 | 0.6032 |
| Anaheim | 0.472 | 0.4578 | 0.4578 | 0.4578 |
| Epcot | 0.5033 | 0.472 | 1 | 0.4578 |
| Duckwater | 0.5802 | 0.539 | 0.6032 | 0.4872 |
| Santa Cruz | 0.5033 | 0.5033 | 0.4578 | 0.4578 |

**Table V.** Results of grey relational grades obtained from models (11) and (13) and the corresponding compromise grades[a]

| Site | $\Gamma_k$ | $\Gamma_k'$ | $\Delta_k(\beta = 0.5)$ |
|---|---|---|---|
| Nome | 0.7515 | 1.1554 | 0.3848 (8) |
| Newark | 1.0000 | 1.0000 | 0.5000 (7) |
| Rock Springs | 1.0000 | 1.3618 | 0.9480 (2) |
| Duquesne | 0.8770 | 1.2808 | 0.6954 (5) |
| Gary | 1.0000 | 1.1642 | 0.7033 (3) |
| Yakima | 0.6642 | 1.068 | 0.1684 (10) |
| Turkey | 1.0000 | 1.123 | 0.6523 (6) |
| Wells | 1.0000 | 1.4038 | 1.0000 (1) |
| Anaheim | 0.5962 | 1.0000 | 0.0000 (12) |
| Epcot | 1.0000 | 1.1609 | 0.6992 (4) |
| Duckwater | 0.6960 | 1.0999 | 0.2474 (9) |
| Santa Cruz | 0.6251 | 1.0289 | 0.0716 (11) |

**Note:** [a] The site ranks are given in parentheses

cost and lives lost. Therefore, one of the advantages of the proposed approach is revealing such alternatives.

*Comparing proposed model and Wu-Olson model*

The model of Wu and Olson is similar to the CCR (after Charnes *et al.* 1978) model (after Charnes, Cooper, & Rhodes, 1978) without explicit inputs. Their model can be obtained by setting $e_j = 0$ and $w_0 = 0$ in model (11). By setting $e_j = 0$, each alternative is allowed to choose its own most favorable weights without using a *priori* weighting and by setting $w_0 = 0$, the assumption of *constant-returns-to-scale* (CRS) in the CCR models is satisfied. Nonetheless, if an alternative has the largest grey relational coefficient in comparison to the other alternatives for a certain attribute, this alternative would always obtain a grey relational grade of 1, even if it has an extremely small grey relational coefficient for other attributes (see Appendix for the mathematical proof). This may lead to the situation in which a large number of alternatives are ranked in the first position. To avoid this issue, we propose the corresponding pessimistic formation by setting $e_j = 0$ and $w'_0 = 0$ in model (13) in which each alternative is allowed to choose its own least favorable weight without using a *priori* weighting under CRS. Under these assumptions, the results obtained from models (11) and (13), as well as the corresponding composite grades for nuclear waste dump site selection, are presented in Table VI.

Table VI shows that Rock Springs , with a compromise grade of 1, stands in the first place, while seven other alternatives are ranked in the first position by model (11). This indicates that the proposed approach can significantly improve the degree of discrimination among alternatives. It is worth noting that, although Rock Springs has the highest compromise grade (=1), it does not have the highest grey relational coefficient with respect to each attribute (Table IV). It is likely due to the fact that Rock Springs not only has relatively high values of grey relational coefficients but also has a better combination among the different attributes.

| Site | $\Gamma_k$ | $\Gamma'_k$ | $\Delta_k(\beta = 0.5)$ |
|---|---|---|---|
| Nome | 1.0000 | 1.0000 | 0.5000 (6) |
| Newark | 1.0000 | 1.0000 | 0.5000 (6) |
| Rock Springs | 1.0000 | 1.7294 | 1.0000 (1) |
| Duquesne | 0.8921 | 1.3176 | 0.5912 (3) |
| Gary | 1.0000 | 1.1146 | 0.5785 (4) |
| Yakima | 0.7855 | 1.0642 | 0.2926 (7) |
| Turkey | 1.0000 | 1.0642 | 0.5440 (5) |
| Wells | 1.0000 | 1.3176 | 0.7177 (2) |
| Anaheim | 0.5735 | 1.0000 | 0.0000 (10) |
| Epcot | 1.0000 | 1.0000 | 0.5000 (6) |
| Duckwater | 0.7351 | 1.0642 | 0.2335 (8) |
| Santa Cruz | 0.5943 | 1.0000 | 0.0244 (9) |

**Note:** [a] The site ranks are given in parentheses

**Table VI.** Results of grey relational grades obtained from models (11) and (13) and the corresponding compromise grades without using a *priori* weighting under CRS[a]

## Conclusions

In this paper, we present the integration of AHP and DEA models in a fuzzy GRA methodology to obtain the weights of attributes. The two sets of grey relational grades obtained from the most and least favorable weights measure the two extreme performances of each alternative. However, any assessment approach considering only one of them is biased. The compromise grey relational grade integrates both the optimistic and the pessimistic performances of each alternative and is, therefore, more comprehensive than either of them. An illustrative example of a nuclear waste dump site selection shows that the compromised grade has a better discriminating power than the optimistic and pessimistic DEA models. We point out that the DEA models discussed in this paper are all based on the so-called weighted additive DEA models, which can be represented in the envelopment forms or the dual multiplier forms. In the envelopment forms, the priority weights of attributes are attached to the slack variables of the objective functions and in the multiplier forms, the upper weight bounds are imposed on the attribute weights.

Therefore, choosing a *priori* weights of attributes, using AHP, in the proposed models is an important matter. However, one of the problems that may occur in practical situations is the difficulty of gathering the different views from some experts for use in AHP. This restricts us in deriving the priority weights of attributes, which further should be used in the proposed additive DEA-based GRA models. At the same time, the equal weight assumption might not be acceptable for decision-makers. In such situations, variable (data-dependent) weights are recommended like range-adjusted weights set introduced in Cooper *et al.* (1999) or a slacks-based distance function proposed by Tone (2001). Finally, the further studies may use those combined AHP and DEA methodologies in conjunction with GRA, which do not necessarily impose weight bounds on the attribute weights. The interested readers may refer to the following papers: Ho and Oh (2010), Jablonsky (2007), Sinuany-Stern *et al.* (2000) for ranking the efficient/inefficient units in DEA models using AHP in a two-stage process. Chen (2002), Cai and Wu (2001), Feng *et al.* (2004), Kim (2000), Pakkar (2014) for weighting the inputs and outputs in the DEA structure using AHP and (Liu and Chen, 2004) for constructing a convex combination of weights using AHP and DEA.

## References

Ahmad, N., Berg, D. and Simons, G.R. (2006), "The integration of analytical hierarchy process and data envelopment analysis in a multi-criteria decision-making problem", *International Journal of Information Technology & Decision Making*, Vol. 5 No. 2, pp. 263-276.

Andersen, P. and Petersen, N.C. (1993), "A procedure for ranking efficient units in data envelopment analysis", *Management Science*, Vol. 39 No. 10, pp. 1261-1264.

Birgün, S. and Güngör, C. (2014), "A multi-criteria call center site selection by hierarchy grey relational analysis", *Journal of Aeronautics and Space Technologies*, Vol. 7 No. 1, pp. 45-52.

Bruce Ho, C. (2011), "Measuring dot com efficiency using a combined DEA and GRA approach", *The Journal of the Operational Research Society*, Vol. 62 No. 4, pp. 776-783.

Cai, Y. and Wu, W. (2001), "Synthetic financial evaluation by a method of combining DEA with AHP", *International Transactions in Operational Research*, Vol. 8 No. 5, pp. 603-609.

Charnes, A., Cooper, W.W. and Rhodes, E. (1978), "Measuring the efficiency of decision making units", *European Journal of Operational Research*, Vol. 2 No. 6, pp. 429-444.

Charnes, A., Cooper, W.W., Golany, B., Seiford, L. and Stutz, J. (1985), "Foundations of data envelopment analysis for Pareto-Koopmans efficient empirical production functions", *Journal of Econometrics*, Vol. 30 No. 1, pp. 91-107.

Chen, T.Y. (2002), "Measuring firm performance with DEA and prior information in Taiwan's banks", *Applied Economics Letters*, Vol. 9 No. 3, pp. 201-204.

Cooper, W.W., Park, K.S. and Pastor, J.T. (1999), "RAM: a range adjusted measure of inefficiency for use with additive models, and relations to other models and measures in DEA", *Journal of Productivity analysis*, Vol. 11 No. 1, pp. 5-42.

Cooper, W.W., Seiford, L.M. and Zhu, J. (Eds) ( 2011), *Handbook on Data Envelopment Analysis*, Springer Science & Business Media, Vol. 164.

Deng, J.L. (1982), "Control problems of grey systems", *Systems & Control Letters*, Vol. 1 No. 5, pp. 288-294.

Dyer, J.S. (1990), "Remarks on the analytic hierarchy process", *Management Science*, Vol. 36 No. 3, pp. 249-258.

Entani, T., Ichihashi, H. and Tanaka, H. (2004), "Evaluation method based on interval AHP and DEA", *Central European Journal of Operations Research*, Vol. 12 No. 1, pp. 25-34.

Feng, Y., Lu, H. and Bi, K. (2004), "An AHP/DEA method for measurement of the efficiency of R&D management activities in universities", *International Transactions in Operational Research*, Vol. 11 No. 2, pp. 181-191.

Girginer, N., Köse, T. and Uçkun, N. (2015), "Efficiency analysis of surgical services by combined use of data envelopment analysis and gray relational analysis", *Journal of Medical Systems*, Vol. 39 No. 5, pp. 1-9.

Goyal, S. and Grover, S. (2012), "Applying fuzzy grey relational analysis for ranking the advanced manufacturing systems", *Grey Systems: Theory and Application*, Vol. 2 No. 2, pp. 284-298.

Hatami-Marbini, A., Saati, S. and Tavana, M. (2013), "Data envelopment analysis with fuzzy parameters: an interactive approach", *Optimizing, Innovating, and Capitalizing on Information Systems for Operations*, pp. 94-108.

Ho, C.B. and Oh, K.B. (2010), "Selecting internet company stocks using a combined DEA and AHP approach", *International Journal of Systems Science*, Vol. 41 No. 3, pp. 325-336.

Hou, J. (2010), "Grey relational analysis method for multiple attribute decision making in intuitionistic fuzzy setting", *Journal of Convergence Information Technology*, Vol. 5 No. 10, pp. 194-199.

Jablonsky, J. (2007), "Measuring the efficiency of production units by AHP models", *Mathematical & Computer Modelling*, Vol. 46 No. 7, pp. 1091-1098.

Jahanshahloo, G.R. and Afzalinejad, M. (2006), "A ranking method based on a full-inefficient frontier", *Applied Mathematical Modelling*, Vol. 30 No. 3, pp. 248-260.

Javanbarg, M.B., Scawthorn, C., Kiyono, J. and Shahbodaghkhan, B. (2012), "Fuzzy AHP-based multicriteria decision making systems using particle swarm optimization", *Expert Systems with Applications*, Vol. 39 No. 1, pp. 960-966.

Jia, W., Li, C. and Wu, X. (2011), "Application of multi-hierarchy grey relational analysis to evaluating natural gas pipeline operation schemes", *Advanced Research on Computer Science and Information Engineering*, Springer Berlin Heidelberg, pp. 245-251.

Jun, G. and Xiaofei, C. (2013), "A coordination research on urban ecosystem in Beijing with weighted grey correlation analysis based on DEA", *Journal of Applied Sciences*, Vol. 13 No. 24, pp. 5749-5759.

Kahraman, C., Cebeci, U. and Ulukan, Z. (2003), "Multi-criteria supplier selection using fuzzy AHP", *Logistics Information Management*, Vol. 16 No. 6, pp. 382-394.

Kim, T. (2000), "Extended topics in the integration of data envelopment analysis and the analytic hierarchy process in decision making", *PhD thesis*, Agricultural & Mechanical College, LA State University.

Kong, W. and Fu, T. (2012), "Assessing the performance of business colleges in Taiwan using data envelopment analysis and student based value-added performance indicators", *Omega*, Vol. 40 No. 5, pp. 541-549.

Kuo, Y., Yang, T. and Huang, G.W. (2008), "The use of grey relational analysis in solving multiple attribute decision-making problems", *Computers & Industrial Engineering*, Vol. 55 No. 1, pp. 80-93.

Lee, A.H.I., Lin, C.Y., Kang, H.Y. and Lee, W.H. (2012), "An integrated performance evaluation model for the photovoltaics industry", *Energies*, Vol. 5 No. 4, pp. 1271-1291.

Lertworasirikul, S., Fang, S.C., Joines, J.A. and Nuttle, H.L. (2003), "Fuzzy data envelopment analysis (DEA): a possibility approach", *Fuzzy Sets and Systems*, Vol. 139 No. 2, pp. 379-394.

Li, G.D., Yamaguchi, D. and Nagai, M. (2008), "A grey-based rough decision-making approach to supplier selection", *The International Journal of Advanced Manufacturing Technology*, Vol. 36 Nos 9/10, pp. 1032-1040.

Liu, C.C. (2009), "A study of optimal weights restriction in data envelopment analysis", *Applied Economics*, Vol. 41 No. 14, pp. 1785-1790.

Liu, C. and Chen, C. (2004), "Incorporating value judgments into data envelopment analysis to improve decision quality for organization", *Journal of American Academy of Business*, Vol. 5 Nos 1/2, pp. 423-427.

Liu, C.M., Hsu, H.S., Wang, S.T. and Lee, H.K. (2005), "A performance evaluation model based on AHP and DEA", *Journal of the Chinese Institute of Industrial Engineers*, Vol. 22 No. 3, pp. 243-251.

Liu, F.H.F. and Chen, C.L. (2009), "The worst-practice DEA model with slack-based measurement", *Computers & Industrial Engineering*, Vol. 57 No. 2, pp. 496-505.

Liu, W.B., Zhang, D.Q., Meng, W., Li, X.X. and Xu, F. (2011), "A study of DEA models without explicit inputs", *Omega*, Vol. 39 No. 5, pp. 472-480.

Markabi, M.S. and Sabbagh, M. (2014), "A hybrid method of grey relational analysis and data envelopment analysis for evaluating and selecting efficient suppliers plus a novel ranking method for grey numbers", *Journal of Industrial Engineering and Management*, Vol. 7 No. 5, pp. 1197-1221.

Markabi, M.S. and Sarbijan, M.S. (2015), "A hybrid model of grey relational analysis and DEA cross-efficiency for the evaluation of decision making units", *International Journal of Economy, Management and Social Science*, Vol. 4 No. 3, pp. 317-322.

Olson, D.L. and Wu, D. (2006), "Simulation of fuzzy multiattribute models for grey relationships", *European Journal of Operational Research*, Vol. 175 No. 1, pp. 111-120.

Osman, M.S., Gadalla, M.H. and Rabie, R.M. (2013), "Fuzzy analytic hierarchical process to determine the relative weights in multi-level programming problems", *International Journal of Mathematical Archive*, Vol. 4 No. 7, pp. 282-295.

Pakkar, M.S. (2014), "Using the AHP and DEA methodologies for stock selection", in Charles, V. and Kumar, M. (Eds), *Business Performance Measurement and Management*, Cambridge Scholars Publishing, Newcastle upon Tyne, pp. 566-580.

Pakkar, M.S. (2016a), "Multiple attribute grey relational analysis using DEA and AHP", *Complex & Intelligent Systems*, (in press), doi: 10.1007/s40747-016-0026-4.

Pakkar, M.S. (2016b), "Using DEA and AHP for hierarchical structures of data", *Industrial Engineering and Management Systems*, Vol. 15 No. 1, pp. 49-62.

Podinovski, V.V. (2004), "Suitability and redundancy of non-homogeneous weight restrictions for measuring the relative efficiency in DEA", *European Journal of Operational Research*, Vol. 154 No. 2, pp. 380-395.

Premachandra, I.M. (2001), "Controlling factor weights in data envelopment analysis by incorporating decision maker's value judgement: an approach based on AHP", *Journal of Information and Management Science*, Vol. 12 No. 2, pp. 1-12.

Saaty, R.W. (1987), "The analytic hierarchy process – what it is and how it is used", *Mathematical Modelling*, Vol. 9 No. 3, pp. 161-176.

Shang, J. and Sueyoshi, T. (1995), "Theory and methodology- a unified framework for the selection of a flexible manufacturing system", *European Journal of Operational Research*, Vol. 85 No. 2, pp. 297-315.

Sinuany-Stern, Z., Mehrez, A. and Hadada, Y. (2000), "An AHP/DEA methodology for ranking decision making units", *International Transactions in Operational Research*, Vol. 7 No. 2, pp. 109-124.

Swim, L.K. (2001), "Improving decision quality in the analytic hierarchy process implementation through knowledge management strategies", *PhD thesis*, The University of Oklahoma.

Takamura, Y. and Tone, K. (2003), "A comparative site evaluation study for relocating Japanese government agencies out of Tokyo", *Socio-Economic Planning Sciences*, Vol. 37 No. 2, pp. 85-102.

Tone, K. (2001), "A slacks-based measure of efficiency in data envelopment analysis", *European Journal of Operational Research*, Vol. 130 No. 3, pp. 498-509.

Tseng, W., Yang, C. and Wang, D. (2009), "Using the DEA and AHP methods on the optimal selection of IT strategic alliance partner", *Proceedings of the 2009 International Conference on Business and Information (BAI 2009)*, Academy of Taiwan Information Systems Research (ATISR), *Kuala Lumpur*, pp. 1-15.

Wang, R.T., Ho, C.T.B. and Oh, K. (2010), "Measuring production and marketing efficiency using grey relation analysis and data envelopment analysis", *International Journal of Production Research*, Vol. 48 No. 1, pp. 183-199.

Wang, Y.M. and Luo, Y. (2006), "DEA efficiency assessment using ideal and anti-ideal decision making units", *Applied Mathematics and Computation*, Vol. 173 No. 2, pp. 902-915.

Wang, Y.M., Chin, K.S. and Yang, J.B. (2007), "Measuring the performances of decision-making units using geometric average efficiency", *Journal of the Operational Research Society*, Vol. 58 No. 7, pp. 929-937.

Wei, G.W. (2010), "GRA method for multiple attribute decision making with incomplete weight information in intuitionistic fuzzy setting", *Knowledge-Based Systems*, Vol. 23 No. 3, pp. 243-247.

Wei, G., Wang, H.J., Lin, R. and Zhao, X. (2011), "Grey relational analysis method for intuitionistic fuzzy multiple attribute decision making with preference information on alternatives", *International Journal of Computational Intelligence Systems*, Vol. 4 No. 2, pp. 164-173.

Wen, M. and Li, H. (2009), "Fuzzy data envelopment analysis (DEA): model and ranking method", *Journal of Computational and Applied Mathematics*, Vol. 223 No. 2, pp. 872-878.

Wu, D.D. and Olson, D.L. (2010), "Fuzzy multiattribute grey related analysis using DEA", *Computers & Mathematics with Applications*, Vol. 60 No. 1, pp. 166-174.

Yang, Y. and John, R. (2012), "Grey sets and greyness", *Information Sciences*, Vol. 185 No. 1, pp. 249-264.

Zeng, G., Jiang, R., Huang, G., Xu, M. and Li, J. (2007), "Optimization of wastewater treatment alternative selection by hierarchy grey relational analysis", *Journal of Environmental Management*, Vol. 82 No. 2, pp. 250-259.

Zhang, J., Wu, D. and Olson, D.L. (2005), "The method of grey related analysis to multiple attribute decision making problems with interval numbers", *Mathematical and Computer Modelling*, Vol. 42 No. 9, pp. 991-998.

Zheng, X. and Lianguang, M. (2013), "Analysis method and its application of weighted grey relevance based on super efficient DEA", *Research Journal of Applied Sciences, Engineering and Technology*, Vol. 5 No. 2, pp. 470-474.

Zhou, P., Ang, B.W. and Poh, K.L. (2007), "A mathematical programming approach to constructing composite indicators", *Ecological Economics*, Vol. 62 No. 2, pp. 291-297.

## Appendix

We assume that $e_j = 0$ and $w_0 = 0$ in model (11). In addition, without loss of generality, we assume that the first alternative has the highest grey relational coefficient within a group of alternatives for the first attribute, that is $\xi_{11} = \max\{\xi_{i1}, i = 1, 2, \ldots, m\}$. Obviously $w_1 = 1/\zeta_{11}, w_2 = \ldots = w_n = 0$ is a feasible solution to model (11) for the first alternative. As $w_1\zeta_{11} + w_2\zeta_{12} + \ldots + w_n\zeta_{1n} - w_0 = 1$ the set of weights is also an optimal solution to model (11) for the first alternative. If the set of weights is used, the first alternative will always obtain a grey relational grade of 1.

**Corresponding author**
Mohammad Sadegh Pakkar can be contacted at: ms_pakkar@laurentian.ca